\documentclass[11pt,reqno]{amsart}
\usepackage{amssymb,amscd,amsbsy}
\usepackage{amssymb,amscd,amsbsy,mathrsfs}
\renewcommand{\a}{\alpha}

\renewcommand{\d}{\delta}
\newcommand{\e}{\varepsilon}
\newcommand{\vk}{\varkappa}

\renewcommand{\l}{\lambda}

\newcommand{\s}{\sigma}

\newcommand{\f}{\varphi}
\renewcommand{\o}{\omega}

\newcommand{\D}{\Delta}

\renewcommand{\O}{\Omega}

\newcommand{\F}{{\mathscr F}}

\newcommand{\h}{{\mathscr H}}

\newcommand{\cL}{{\mathscr L}}
\newcommand{\M}{{\mathscr M}}

\newcommand{\X}{{\mathscr X}}
\newcommand{\Y}{{\mathscr Y}}

\newcommand{\R}{{\Bbb R}}
\newcommand{\Z}{{\Bbb Z}}

\newcommand{\0}{{\boldsymbol{0}}}

\newcommand{\bs}{\boldsymbol}

\newcommand{\m}{{\boldsymbol m}}
\newcommand{\bS}{{\boldsymbol S}}

\newcommand{\rf}[1]{(\ref{#1})}

\newcommand{\df}{\stackrel{\mathrm{def}}{=}}

\newcommand{\supp}{\operatorname{supp}}

\newcommand{\trace}{\operatorname{trace}}

\newcommand{\const}{\operatorname{const}}
\newcommand{\tr}{\operatorname{trace}}
\newcommand{\eeq}{\end{equation}}
\newcommand{\beq}{\begin{equation}}
\newcommand{\bay}{\begin{eqnarray}}
\newcommand{\ba}{\begin{align*}}
\newcommand{\ea}{\end{align*}}
\newcommand{\ey}{\end{eqnarray}}
\newcommand{\bey}{\begin{eqnarray*}}
\newcommand{\eey}{\end{eqnarray*}}

\newcommand{\be}{\infty}

\newcommand{\bl}{\blacksquare}

\newcommand{\Pf}{{\bf Proof. }}

\newtheorem{thm}{\hspace{\parindent}Theorem}[section]

\newtheorem{cor}[thm]{\hspace{\parindent}Corollary}
\newtheorem{lem}[thm]{\hspace{\parindent}Lemma}

\theoremstyle{remark}

\newtheorem*{rem*}{Remark}

\newcommand\cZ{\mathcal{Z}}
\newcommand\dg{\frak D}

\newcommand{\Tp}{{\mathscr T}^{(m)}_{A,K}}

\newcommand{\sm}{{\mathcal S}}

\begin{document}

\newcommand{\vse}{\vspace{.2in}}
\numberwithin{equation}{section}

\title{Trace formulae for perturbations of class $\bs{\bS_m}$}
\author{A.B. Aleksandrov and V.V. Peller}
\thanks{The first author is partially supported by RFBR grant 08-01-00358-a and by
Russian Federation presidential grant NSh-2409.2008.1;
the second author is partially supported by NSF grant DMS 1001844 and by ARC grant}

\newcommand{\mt}{{\mathcal T}}

\begin{abstract}
We obtain general trace formulae in the case of perturbation of self-adjoint operators by self-adjoint operators of class $\bS_m$, where $m$ is a positive integer. In \cite{PSS} a trace formula for operator Taylor polynomials was obtained. This formula includes the Livshits--Krein trace formula in the case $m=1$ and the Koplienko trace formula in the case $m=2$.
We establish most general trace formulae in the case of perturbation of Schatten--von Neumann class $\bS_m$. We also improve the trace formula obtained in \cite{PSS} for operator Taylor polynomials and prove it for arbitrary functions in he Besov space $B_{\be1}^m(\R)$.

We consider several other special cases of our general trace formulae. In particular, we establish a trace formula for $m$th order operator differences.
\end{abstract}

\maketitle

\

\setcounter{section}{0}
\section{\bf Introduction}
\setcounter{equation}{0}

\

The purpose of this paper is to obtain most general trace formulae for perturbation of self-adjoint operators by operators of class $\bS_m$, where $m$ is a positive integer.

The Livshits--Krein spectral shift function for trace class perturbations of self-adjoint operators was introduced by Livshits \cite{L} in a special case and by M.G. Krein \cite{Kr} in the general case. It was proved in \cite{Kr} that for a
(not necessarily bounded) self-adjoint operator $A$ and a trace class self-adjoint operator $K$
there exists a unique function $\xi$ in $L^1(\R)$ such that
\bay
\label{LKf}
\trace\big(f(A+K)-f(A)\big)=\int_\R f'(x)\xi(x)\,dx
\ey
for every function $f$ such that the Fourier transform $\F f'$ of its derivative belongs to $L^1(\R)$. Note that the right-hand side of \rf{LKf} is well defined
for every Lipschitz functions $f$. Krein conjectured that for every Lipschitz function $f$ the operator $f(A+K)-f(A)$ belongs to the trace class $\bS_1$ and \rf{LKf} holds. It turns out that this is wrong. Farforovskaya constructed a counter-example in \cite{F1}.

Later in \cite{Pe1} a necessary condition was found. Namely, it was shown in
\cite{Pe1} that if $f(A+K)-f(A)\in\bS_1$ for every self-adjoint $A$ and every self-adjoint $K$ in $\bS_1$, then $f$ locally belongs to the Besov space $B_{11}^1(\R)$
(see \S\,2 for a brief introduction in Besov spaces).
This necessary condition also implies that there are Lipschitz functions $f$, for which the condition $K\in\bS_1$ does not imply that $f(A+K)-f(A)\in\bS_1$.

On the other hand, it was shown in \cite{Pe1} and \cite{Pe2} that if $f$ belongs to the Besov class $B_{\be1}^1(\R)$, then the left-hand side of \rf{LKf} belongs to $\bS_1$ and trace formula \rf{LKf} holds.

We refer the reader to the survey article \cite{BY} for more information about the Livshits--Krein trace formula.

Koplienko considered in \cite{Ko} the case of perturbations by self-adjoint operators of Hilbert--Schimidt class $\bS_2$. With each pair of self-adjoint operators $(A,K)$ such that $K\in\bS_2$ he associated a unique nonnegative function $\eta$ in $L^1(\R)$ such that
\bay
\label{Ktf}
\trace\left(f(A+K)-f(A)-\frac{d}{dt}f(A+tK)\Big|_{t=0}\right)
=\int_\R f''(x)\eta(x)\,dx
\ey
for every rational function $f$ bounded on $\R$.

In \cite{Pe4} the result of Koplienko was improved. It was shown in \cite{Pe4} that if $f$ belongs to the Besov space $B_{\be1}^2(\R)$, then the operator on the left-hand side of \rf{Ktf} belongs to $\bS_1$ and formula \rf{Ktf} holds.

Koplienko also attempted in \cite{Ko} to generalize his results to the case of perturbations of Schatten--von Neumann class $\bS_m$ for an arbitrary positive integer $m$. However, his proof for $m>2$ was erroneous.

We also mention here the paper \cite{GPS}, in which interesting results related to the Koplienko trace formula were obtained.

In \cite{PSS} it was shown that for every positive integer $m$
and for every pair $(A,K)$ of self-adjoint operators with $K\in\bS_m$, there exists a unique function
$\eta_m$ in $L^1(\R)$ such that the following trace formula holds:
\bay
\label{PSS}
\trace\left({\mathscr T}^{(m)}_{A,K}f\right)=\int_\R f^{(m)}(x)\eta_m(x)\,dx
\ey
for functions $f$ satisfying the conditions
\bay
\label{PrFo}
\F f^{(j)}\in L^1(\R),\quad 0\le j\le m,
\ey
where the operator Taylor polynomial ${\mathscr T}^{(m)}_{A,K}f$ is defined by
\begin{align*}
{\mathscr T}^{(m)}_{A,K}f&\df f(A+K)-f(A)\\[.2cm]
&-\frac{d}{dt}f(A+tK)\Big|_{t=0}-
\cdots-\frac1{(m-1)!}\frac{d^{m-1}}{dt^{m-1}}f(A+tK)\Big|_{t=0}.
\end{align*}
The function $\eta_m$ is called the {\it spectral shift function of order $m$}.

Note that earlier partial results were obtained in \cite{DS}, \cite{Sk1}, and \cite{Sk2}.

In \S\,6 of this paper we obtain  most general trace formulae that include the trace formula for operator Taylor polynomials as a special case.

Moreover,
we improve the result of \cite{PSS} for operator Taylor polynomials. We prove in \S\,7 that trace formula
\rf{PSS} holds under much less restrictive assumptions on $f$: we show that it holds for all functions $f$ in the Besov space $B_{\be1}^m(\R)$. At the same time we establish our general trace formulae also for arbitrary functions in the Besov space $B_{\be1}^m(\R)$.

In \S\,5 we establish a formula that expresses the operator Taylor polynomial ${\mathscr T}^{(m)}_{A,K}f$ in terms of a multiple operator integral. This formula will not be used to obtain trace formulae. However, we believe that it is of independent interest.


In \S\,2 we give a brief introduction to Besov spaces and in \S 3 we introduce multiple operator integrals. Finally, in \S\,4 we prove two theorems on continuous dependence of multiple operator integrals on the corresponding self-adjoint operators that will be used to obtain our main results.

Throughout the paper we use the notation $\m$ for Lebesgue measure on $\R$.

\

\section{\bf Besov classes}
\setcounter{equation}{0}

\

The purpose of this section is to give a brief introduction to the Besov spaces that play an important role in problems of perturbation theory.

Let $0<p,\,q\le\be$ and $s\in\R$. The homogeneous Besov class $B^s_{pq}(\R)$ of functions (or distributions) on $\R$ can be defined in the following way. Let $w$ be an infinitely differentiable function on $\R$ such
that
$$
w\ge0,\quad\supp w\subset\left[\frac12,2\right],\quad\mbox{and} \quad w(x)=1-w\left(\frac x2\right)\quad\mbox{for}\quad x\in[1,2].
$$
We
define the functions $W_n$ and $W^\sharp_n$ on $\R$ by
$$
\F W_n(x)=w\left(\frac{x}{2^n}\right),\quad\F W^\sharp_n(x)=\F W_n(-x),\quad n\in\Z,
$$
where $\F$ is the {\it Fourier transform}:
$$
\big(\F f\big)(t)=\int_\R f(x)e^{-{\rm i}xt}\,dx,\quad f\in L^1.
$$

With every tempered distribution $f\in{\mathscr S}^\prime(\R)$ we
associate a sequences $\{f_n\}_{n\in\Z}$,
\bay
\label{fn}
f_n\df f*W_n+f*W_n^\sharp.
\ey
Initially we define the (homogeneous) Besov class $\dot B^s_{pq}(\R)$ as the set of all $f\in{\mathscr S}^\prime(\R)$
such that
\bay
\label{Wn}
\{2^{ns}\|f_n\|_{L^p}\}_{n\in\Z}\in\ell^q(\Z).
\ey
According to this definition, the space $\dot B^s_{pq}(\R)$ contains all polynomials. Moreover, the distribution $f$ is defined by the sequence $\{f_n\}_{n\in\Z}$
uniquely up to a polynomial. It is easy to see that the series $\sum_{n\ge0}f_n$ converges in ${\mathscr S}^\prime(\R)$.
However, the series $\sum_{n<0}f_n$ can diverge in general. It is easy to prove that the
series $\sum_{n<0}f_n^{(r)}$ converges on uniformly $\R$ for each nonnegative integer
$r>s-1/p$ if $q>1$ and the series $\sum_{n<0}f_n^{(r)}$ converges uniformly, whenever $r\ge s-1/p$ if $q\le1$.

Now we can define the modified (homogeneous) Besov class $B^s_{pq}(\R)$. We say that a distribution $f$
belongs to $B^s_{pq}(\R)$ if $\{2^{ns}\|f_n\|_{L^p}\}_{n\in\Z}\in\ell^q(\Z)$ and
$f^{(r)}=\sum_{n\in\Z}f_n^{(r)}$ in the space ${\mathscr S}^\prime(\R)$, where
$r$ is the minimal nonnegative integer such that $r>s-1/p$ in the case $q>1$ and
$r$ is the minimal nonnegative integer such that $r\ge s-1/p$ in the case $q\le1$.
Now the function $f$ is determined uniquely by the sequence $\{f_n\}_{n\in\Z}$ up
to a polynomial of degree less that $r$, and a polynomial $\varphi$ belongs to $B^s_{pq}(\R)$
if and only if $\deg\varphi<r$.



Besov spaces $B^s_{pq}(\R)$ admit equivalent definitions in terms of finite
differences. We give such a definition in the case of Besov spaces
$B_{\be1}^m(\R)$, with which we mostly deal in this paper.

For $t\in\R$, we define the difference operator $\D_t$ by
$$
\big(\D_tf\big)(x)\df f(x+t)-f(x).
$$
A function $f$ belongs to $B_{\be1}^m(\R)$ if and only if
$$
\int_\R\frac{\|\D_t^{m+1}f\|_{L^\be}}{|t|^{1+m}}\,dt<\be\quad\mbox{and}\quad
\lim_{|x|\to\be}\frac{|f(x)|}{(1+|x|)^m}=0.
$$

We refer the reader to \cite{Pee} and \cite{Pe3} for more detailed information on Besov spaces.

\

\section{\bf Multiple operator integrals}
\setcounter{equation}{0}
\label{koi}

\

In this section we give a brief introduction to the theory of multiple operator integrals. Double operator integrals appeared in the paper \cite{DK} by Daletskii and S.G. Krein. However, the beautiful theory of double operator integrals was developed later by Birman and Solomyak in \cite{BS1}, \cite{BS2}, and \cite{BS3}.

We are not going to define double operator integrals
$\iint\Phi(x,y)\,d E_1(x)\,T\,dE_2(y)$ in the case of Hilbert--Schmidt
operators $T$ that was the starting point in \cite{BS1}, \cite{BS2}, and \cite{BS3}. We use the approach based on (integral) projective tensor products.
In the case of bounded or trace class operators $T$ this approach is equivalent to the approach of Birman and Solomyak, see \cite{Pe1}. Moreover, we start with the definition of the more general notion of multiple operator integrals. Multiple operator integrals were defined in terms of integral projective tensor products in \cite{Pe5} (see earlier publications \cite{Pa} and \cite{St}, where multiple operators were defined under much more restrictive assumptions).

To simplify the notation, we consider the case of triple operator integrals. The definition for general multiple operator integrals is the same.
Let $(\X,E_1)$, $(\Y,E_2)$, and $(\cZ,E_3)$ be spaces with spectral measures $E_1$ and $E_2$, and $E_3$
on Hilbert spaces $\h_1$ and $\h_2$, and $\h_3$. Suppose that a function $\Phi$ on $\X\times\Y\times\cZ$ belongs to the {\it projective tensor
product}
$L^\be(E_1)\hat\otimes L^\be(E_2)\hat\otimes L^\be(E_3)$ of $L^\be(E_1)$, $L^\be(E_2)$, and
$L^\be(E_3)$ (i.e., $\Phi$ admits a representation
\bay
\label{ptp}
\Phi(x,y)=\sum_{n\ge0}\f_n(x)\psi_n(y)\chi_n(z),
\ey
where $\f_n\in L^\be(E_1)$, $\psi_n\in L^\be(E_2)$, and $\chi_n\in L^\be(E_3)$ are
functions such that
\bay
\label{ptpn}
\sum_{n\ge0}\|\f_n\|_{L^\be}\|\psi_n\|_{L^\be}\|\chi_n\|_{L^\be}<\be).
\ey
Then for arbitrary bounded linear operators $T_1:\h_2\to\h_1$ and $T_2:\h_3\to\h_2$, we put
\begin{align*}
\int\limits_\X\int\limits_\Y\int\limits_\cZ
&\Phi(x,y,z)\,d E_1(x)T_1\,dE_2(y)T_2\,dE_3(z)
\\[.2cm]
&\df
\sum_{n\ge0}\left(\,\int\limits_\X \f_n\,dE_1\right)T_1\left(\,\int\limits_\Y \psi_n\,dE_2\right)T_2\left(\int\limits_\cZ\chi_n\,dE_3\right).
\end{align*}
It was shown in \cite{Pe5} (see also \cite{ACDS} for a different proof)  that the above definition does not depend on the choice of a representation \rf{ptp}.
For $\Phi\in L^\be(E_1)\hat\otimes L^\be(E_2)\hat\otimes L^\be(E_3)$, its norm
is, by definition, the infimum of the  left-hand side of \rf{ptpn} over all representations \rf{ptp}.
\newcommand{\fI}{\mathcal I}

We can enlarge the class of functions $\Phi$, for which multiple operator integrals can be defined by considering {\it integral projective tensor products}.
This approach for multiple operator integrals was given in \cite{Pe5}.
Again, we consider here the case of triple operator integrals; the general case can be treated in the same way.

We say that a measurable function $\Phi$ on $\X\times\Y\times\cZ$
belongs to the {\it integral projective tensor product}
$L^\be(E_1)\hat\otimes_{\rm i}L^\be(E_2)\hat\otimes_{\rm i}L^\be(E_3)$ if $\Phi$ admits a representation
\bay
\label{ttp}
\Phi(x,y,z)=\int_\O \f(x,\o)\psi(y,\o)\chi(z,\o)\,d\vk(\o),
\ey
where $(\O,\vk)$ is a measure space with a $\s$-finite measure, $\f$ is a measurable function on $\X\times \O$,
$\psi$ is a measurable function on $\Y\times \O$, $\chi$ is a measurable function on $\cZ\times \O$,
and
\bay
\label{iptpn}
\int_\O\|\f(\cdot,\o)\|_{L^\be(E_1)}\|\psi(\cdot,\o)\|_{L^\be(E_2)}\|\chi(\cdot,\o)\|_{L^\be(E_3)}\,d\vk(\o)<\be.
\ey

For a bounded linear operator $T_1$ from $\h_2$ to $\h_1$, a bounded linear operator $T_2$ from $\h_3$ to $\h_2$, and a function $\Phi$ in
$L^\be(E_1)\hat\otimes_{\rm i}L^\be(E_2)\hat\otimes_{\rm i}L^\be(E_3)$ of the form \rf{ttp}, we put
\begin{align}
\label{opr}
&\int\limits_\X\int\limits_\Y\int\limits_\cZ\Phi(x,y,z)
\,d E_1(x)T_1\,dE_2(y)T_2\,dE_3(z)\\[.2cm]
\df&\int\limits_\O\left(\,\int\limits_\X \f(x,\o)\,dE_1(x)\right)T_1
\left(\,\int\limits_\Y \psi(y,\o)\,dE_2(y)\right)T_2
\left(\,\int\limits_\cZ \chi(z,\o)\,dE_3(z)\right)\,d\vk(\o).\nonumber
\end{align}

Again, the above definition does not depend on the choice of a representation \rf{ttp} (see \cite{Pe5}). The norm $\|\Phi\|_{L^\be\hat\otimes_{\rm i}L^\be\hat\otimes_{\rm i}L^\be}$ is defined as the infimum of the left-hand side of \rf{iptpn} over all representations \rf{ttp}.

It is easy to see that the following inequality holds
$$
\left\|\int\limits_\X\int\limits_\Y\int\limits_\cZ\Phi(x,y,z)
\,dE_1(x)T_1\,dE_2(y)T_2\,dE_3(z)\right\|
\le\|\Phi\|_{L^\be\hat\otimes_{\rm i}L^\be\hat\otimes_{\rm i}L^\be}\|T_1\|\cdot\|T_2\|.
$$

It is also easy to see that if $T_1\in\bS_p$ and $T_2\in\bS_q$, and $1/p+1/q\le1$,
then the triple operator integral \rf{opr} belongs to $\bS_r$ and
$$
\left\|\int\limits_\X\int\limits_\Y\int\limits_\cZ\Phi(x,y,z)
\,dE_1(x)T_1\,dE_2(y)T_2\,dE_3(z)\right\|_{\bS_r}
\le\|\Phi\|_{L^\be\hat\otimes_{\rm i}L^\be\hat\otimes_{\rm i}L^\be}
\|T_1\|_{\bS_p}\cdot\|T_2\|_{\bS_q},
$$
where $1/r=1/p+1/q$.

Note that in \cite{JTT} Haagerup tensor products were used to define multiple operator integrals. However, it is not clear whether this can lead to stronger results in perturbation theory.

For a function $f$ on $\R$ the  {\it divided differences $\dg^k f$ of order $k$} are defined inductively as follows:
$$
\dg^0f\df f;
$$
if $k\ge1$, then in the case when $t_1,t_2,\cdots,t_{k+1}$ are distinct points of $\R$,
$$
(\dg^{k}f)(x_1,\cdots,x_{k+1})\df
\frac{(\dg^{k-1}f)(x_1,\cdots,x_{k-1},x_k)-
(\dg^{k-1}f)(x_1,\cdots,x_{k-1},x_{k+1})}{x_{k}-x_{k+1}}
$$
(the definition does not depend on the order of the variables). We put
$$
\dg f=\dg^1f.
$$
If $f\in C^k(\R)$, then $\dg^{k}f$ extends by continuity to a function defined for all points $x_1,x_2,\cdots,x_{k+1}$.

It can be shown that
$$
({\frak D}^m f)(x_1,\dots,x_{m+1})=\sum\limits_{k=1}^{m+1}f(x_k)
\prod\limits_{j=1}^{k-1}(x_k-x_j)^{-1}\prod\limits_{j=k+1}^{m+1}(x_k-x_j)^{-1}.
$$

It follows from the results of Birman and Solomyak \cite{BS3} that if
 $A$ is a self-adjoint operator (not necessarily bounded),
$K$ is a bounded self-adjoint operator, and
$f$ is a continuously differentiable
function on $\R$ such that
$\dg f\in L^\be(E_{A+K})\hat\otimes_{\rm i} L^\be(E_A)$, then
\bay
\label{BSF}
f(A+K)-f(A)=\iint\limits_{\R\times\R}\big(\dg f\big)(x,y)\,dE_{A+K}(x)K\,dE_A(y)
\ey
(throughout this paper $E_A$ stands for the spectral measure of $A$).


It was shown in \cite{Pe5} that if $f$ is a bounded function on $\R$ such that
its Fourier transform $\F f$ is supported on $[-\s,\s]$, then
then
\bay
\label{Dmf}
\dg^mf\in\underbrace{L^\be(\R)\hat\otimes_{\rm i}\cdots\hat\otimes_{\rm i} L^\be(\R)}_{m+1}
\ey
and
\bay
\label{Boke}
\big\|\dg^m f\big\|_{L^\be\hat\otimes_{\rm i}\cdots\hat\otimes_{\rm i} L^\be}\le\const \s^m\|f\|_{L^\be}.
\ey
This implies (see \cite{Pe5}) that if $f$ belongs to the Besov space $B_{\be1}^m(\R)$,
then \rf{Dmf} holds
and
\bay
\label{Besov}
\big\|\dg^m f\big\|_{L^\be\hat\otimes_{\rm i}\cdots\hat\otimes_{\rm i} L^\be}\le\const \|f\|_{B_{\be1}^m(\R)}.
\ey

It was also proved in \cite{Pe5} that if $f\in B_{\be1}^m(\R)$, $A$ is a self-adjoint operator and $K$ is a bounded self-adjoint operator, then
the function
$$
t\mapsto f(A+tK)
$$
has $m$ derivatives in the operator norm and
\bay
\label{mpro}
\frac{d^m}{dt^m}f(A+tK)\Big|_{t=0}\!
=\!\underbrace{\int\!\!\cdots\!\!\int}_{m+1}\!\big(\dg^mf\big)\big(x_1,\cdots,x_{m+1}\big)
dE_A(x_1)K\cdots KdE_A(x_{m+1}).
\ey
Strictly speaking under the condition $f\in B_{\be1}^m(\R)$ lower derivatives
$$
\frac{d^k}{dt^k}f(A+tK),\quad k<m,
$$
do not necessarily exist. To get the $m$th derivative we should define the left-hand side of \rf{mpro} by
$$
\frac{d^m}{dt^m}f(A+tK)=\sum_{n\in\Z}\,\frac{d^m}{dt^m}f_n(A+tK),
$$
where the functions $f_n$, $n\in\Z$, are defined by \rf{fn}.

\

\section{\bf Continuous dependence of multiple operator integrals}
\setcounter{equation}{0}
\label{vspom}

\

We have already seen (see \rf{mpro}) that multiple operator integrals
that involve divided differences of an arbitrary order play and important role in perturbation theory. In this paper we are going to consider other multiple operator integrals that involve divided differences.

The purpose of this section is to show that such multiple operator integrals depend continuously on the corresponding self-adjoint operators. We establish two similar results. The first one establishes continuous dependence in the operator (or trace class) norm, while the second one deals with continuous dependence in the strong operator topology.

We are going to use a representation of $\dg^m f$ as an element of the integral projective tensor product of $L^\be$ spaces. This representation was obtained in \cite{Pe2} for $m=1$ and in \cite{Pe5} for $m\ge2$. To simplify the notation,
we formulate the result for $m=1$ and $m=2$. Similar formulae hold for $m>2$.

Let $f$ be a bounded function on $\R$ such that
the support of its Fourier transform $\F f$ is a compact subset of $(0,\be)$. Then (see \cite{Pe2})
\begin{align*}
(\dg f)(x,y)&=\frac{\rm i}{2\pi}\iint\limits_{\R_+\times\R_+}(\F f)(s+t)\,e^{{\rm i}sx}e^{{\rm i}ty}\,ds\,dt\\[.2cm]
&=\frac{\rm i}{2\pi}\iint\limits_{\R_+\times\R_+}\frac s{s+t}(\F f)(s+t)\,e^{{\rm i}sx}e^{{\rm i}ty}\,dsdt\\[.2cm]
+&\frac{\rm i}{2\pi}\iint\limits_{\R_+\times\R_+}\frac t{s+t}(\F f)(s+t)\,e^{{\rm i}sx}e^{{\rm i}ty}\,ds\,dt\\[.2cm]
&={\rm i}\int_{\R_+}\big((\sm^{-t}f)*q_t\big)(x)e^{{\rm i}ty}\,dt
+{\rm i}\int_{\R_+}e^{{\rm i}sx}\big((\sm^{-s}f)*q_s\big)(y)\,ds,
\end{align*}
where $(\sm^{-t}f)(x)\df e^{-{\rm i}tx}f(x)$ and $q_t$ is the distribution defined by
$$
(\F q_t)(s)=\max\left\{1-\frac{t}{|s+t|},0\right\}.
%
$$
One can prove that for $t>0$, $q_t$ is a (complex) measure  whose total variation does not depend on $t$, see \cite{Pe2}.

Hence,
\bay
\label{j2}
(\dg f)(x,y)={\rm i}\int_{\R_+}\big((\sm^{-t}f)*q_t\big)(x)e^{{\rm i}ty}\,dt
+{\rm i}\int_{\R_+}e^{{\rm i}sx}\big((\sm^{-s}f)*q_s\big)(y)\,ds
\ey
for every $f\in L^\be(\R)$ such that the support of its Fourier transform $\F f$ is a compact subset of $(0,\be)$. Moreover, if $\supp\F f\subset[0,\s]$, $\s>0$, then
\bay
\label{s1}
\int_{\R_+}\big\|\big((\sm^{-s}f)*q_s\big)\big\|_{L^\be}\,dt\le\const\s\|f\|_{L^\be}.
\ey

In the same way we have

\begin{align*}
(\dg^m f)&(x_1,x_2,\dots,x_{m+1})=
\frac{{\rm i}^m}{2\pi}\int\limits_0^\be\cdots\int\limits_0^\be
(\F f)\left(\sum_{j=1}^{m+1} t_j\right)
e^{{\rm i}t_1x_1}e^{{\rm i}t_2x_2}\dots e^{{\rm i}t_{m+1}x_{m+1}}\,dt\\[.3cm]
=&\frac{{\rm i}^{m}}{2\pi}\sum_{j=1}^{m+1}\,
\int\limits_0^\be\cdots\int\limits_0^\be{t_j}\left(\sum_{j=1}^{m+1} t_j\right)^{-1}
(\F f)\left(\sum_{j=1}^{m+1} t_j\right)e^{{\rm i}t_1x_1}e^{{\rm i}t_2x_2}\dots e^{{\rm i}t_{m+1}x_{m+1}}\,dt
\end{align*}
for every function $f\in C(\R)$ such that $\F f\in L^1(\R)$ and $\supp\F f\subset[0,\be)$.
Applying this formula for $m=2$, we obtain (see \cite{Pe5})
\begin{align}
\label{j3}
(\dg^2f)(x,y,z)=&
-\iint\limits_{\R_+\times\R_+}\big((\sm^{-t-u}f)*q_{t+u}\big)(x)e^{{\rm i}ty}e^{{\rm i}uz}
\,dt\,du\nonumber\\[.2cm]
&-\iint\limits_{\R_+\times\R_+}e^{{\rm i}sx}\big((\sm^{-s-u}f)*q_{s+u}\big)(y)e^{{\rm i}uz}
\,ds\,du\nonumber\\[.2cm]
&-\iint\limits_{\R_+\times\R_+}e^{{\rm i}sx}e^{{\rm i}ty}\big((\sm^{-s-t}f)*q_{s+t}\big)(z)
\,ds\,dt
\end{align}
for every $f\in L^\be(\R)$ such that the support of its Fourier transform $\F f$ is a compact subset of $(0,\be)$.

Moreover, if $\supp\F f\subset[0,\s]$, $\s>0$, then
\bay
\label{s2}
\iint\limits_{\R_+\times\R_+}
\big\|\big((\sm^{-t-u}f)*q_{t+u}\big)\big\|_{L^\be}\,dt\,du
\le\const\s^2\|f\|_{L^\be}.
\ey

\begin{thm}
\label{preper}
Let $m$ be a positive integer and let $f\in B_{\be1}^m(\R)$. Suppose that
$K_1,\,K_2,\cdots,K_m$ are bounded linear operators,
$A_1,\,A_2,\cdots,A_{m+1}$ are (not necessarily bounded) self-adjoint operators, and $\big\{A_{1,j}\big\}_{j\ge0},\,\big\{A_{2,j}\big\}_{j\ge0},\cdots,\big\{A_{m+1,j}\big\}_{j\ge0}$ are sequences of self-adjoint operators such that
$$
\lim_{j\to\be}\|A_1-A_{1,j}\|=\lim_{j\to\be}\|A_2-A_{2,j}\|=\cdots
=\lim_{j\to\be}\|A_{m+1}-A_{m+1,j}\|=0.
$$
Then
\begin{align}
\label{skho}
\!\!\!\!\lim_{j\to\be}&\underbrace{\int\cdots\int}_{m+1}(\dg^{m}f)(x_1,\cdots,x_{m+1})
\,dE_{A_{1,j}}(x_1)K_1\,dE_{A_{2,j}}(x_2)K_2\cdots K_m\,dE_{A_{m+1,j}}(x_{m+1})
\nonumber\\[.2cm]
=&\underbrace{\int\cdots\int}_{m+1}(\dg^{m}f)(x_1,\cdots,x_{m+1})
\,dE_{A_1}(x_1)K_1\,dE_{A_2}(x_2)K_2\cdots K_m\,dE_{A_{m+1}}(x_{m+1})
\end{align}
in the operator norm.

Moreover, if $K_1,\,K_2,\cdots,K_m\in\bS_m$, then the limit in {\em\rf{skho}} exists in the norm of $\bS_1$.
\end{thm}

\Pf The proof is similar to the proof of Theorem 5.1 in \cite{Pe5}.
For simplicity, we give the proof in the case $m=2$. For arbitrary $m$, the proof is the same.

As usual, it suffices to prove the result in the case when
$\supp\F f$ is a compact subset of $(0,\be)$. We are going to use formula
\rf{j3}. Put
$$
\Phi(x_1,x_2,x_3)=
\iint\limits_{\R_+\times\R_+}\big((\sm^{-t-u}f)*q_{t+u}\big)(x_1)e^{{\rm i}tx_2}e^{{\rm i}ux_3}\,dt\,du.
$$

We have to show that
\begin{align*}
\lim_{j\to\be}&\iiint\Phi(x_1,x_2,x_3)
\,dE_{A_{1,j}}(x_1)K_1\,dE_{A_{2,j}}(x_2)K_2\,dE_{A_{3,j}}(x_3)
\nonumber\\[.2cm]
=&\iiint\Phi(x_1,x_2,x_3)
\,dE_{A_1}(x_1)K_1\,dE_{A_2}(x_2)K_2\,dE_{A_3}(x_3).
\end{align*}
The proof for the other two terms in \rf{j3} is exactly the same.

Clearly,
\begin{align*}
\iiint&\Phi(x_1,x_2,x_3)
\,dE_{A_{1,j}}(x_1)K_1\,dE_{A_{2,j}}(x_2)K_2\,dE_{A_{3,j}}(x_3)\\[.2cm]
=\iint\limits_{\R_+\times\R_+}&\big((\sm^{-t-u}f)*q_{t+u}\big)(A_{1,j})K_1e^{{\rm i}tA_{2,j}}K_2e^{{\rm i}uA_{3,j}}\,dt\,du
\end{align*}
and
\begin{align*}
\iiint&\Phi(x_1,x_2,x_3)
\,dE_{A_1}(x_1)K_1\,dE_{A_2}(x_2)K_2\,dE_{A_3}(x_3)\\[.2cm]
=\iint\limits_{\R_+\times\R_+}&\big((\sm^{-t-u}f)*q_{t+u}\big)(A_1)K_1e^{{\rm i}tA_2}K_2e^{{\rm i}uA_3}\,dt\,du.
\end{align*}

It is easy to see that it suffices to show that
$$
\lim_{j\to\be}\big\|\big((\sm^{-t-u}f)*q_{t+u}\big)(A_{1,j})-
\big((\sm^{-t-u}f)*q_{t+u}\big)(A_1)\big\|=0,
$$
$$
\lim_{j\to\be}\big\|e^{{\rm i}tA_{2,j}}-e^{{\rm i}tA_2}\big\|=
\lim_{j\to\be}\big\|e^{{\rm i}uA_{3,j}}-e^{{\rm i}uA_3}\big\|=0
$$
for every $t,\,u>0$. However, this is obvious, since the functions
$(\sm^{-t-u}f)*q_{t+u}$ and $x\mapsto e^{itx}$ are operator Lipschitz.

The same reasoning shows that if $K_1,\,K_2,\cdots,K_m\in\bS_m$, then the limit in \rf{skho} exists in the norm of $\bS_1$.
$\bl$

\begin{thm}
\label{preper0}
Let $m$ be a positive integer and let $f\in B_{\be1}^m(\R)$. Suppose that
$K_1,\,K_2,\dots,K_m$ are bounded linear operators,
$A_1,\,A_2,\dots,A_{m+1}$ are (not necessarily bounded) self-adjoint operators, and $\big\{A_{1,j}\big\}_{j\ge0},\,\big\{A_{2,j}\big\}_{j\ge0},\dots,\big\{A_{m+1,j}\big\}_{j\ge0}$ are sequences of bounded self-adjoint operators such that
$$
\lim_{j\to\be}\|A_kx-A_{k,j}x\|=0
$$
for $k=1,2,\dots,m+1$ and $x$ in the domain of  $A_k$.
Then
\begin{align*}
\!\!\!\!\lim_{j\to\be}&\underbrace{\int\cdots\int}_{m+1}(\dg^{m}f)(x_1,\dots,x_{m+1})
\,dE_{A_{1,j}}(x_1)K_1\,dE_{A_{2,j}}(x_2)K_2\dots K_m\,dE_{A_{m+1,j}}(x_{m+1})
\nonumber\\[.2cm]
=&\underbrace{\int\cdots\int}_{m+1}(\dg^{m}f)(x_1,\dots,x_{m+1})
\,dE_{A_1}(x_1)K_1\,dE_{A_2}(x_2)K_2\dots K_m\,dE_{A_{m+1}}(x_{m+1})
\end{align*}
in the strong operator topology.
\end{thm}

We need two lemmata. The one is Lemma 8.4 in \cite{AP1}.

\begin{lem}
\label{sot}
Let $f$ be a bounded continuous function on $\R$. Suppose that $A$ is a self-adjoint operator (not necessarily bounded) and $\{A_j\}_{j\ge0}$ is a sequence of bounded self-adjoint operators such that
\bey
\label{Aj}
\lim_{j\to\be}\|A_ju-Au\|=0\quad\mbox{for every}\quad u\quad\mbox{in the domain of}\quad A.
\eey
Then
\bey
\label{sil}
\lim_{j\to\be}f(A_j)=f(A)\quad\mbox{in the strong operator topology}.
\eey
\end{lem}

The second lemma is an obvious operator-valued version of
the Lebesgue dominated convergence theorem.

\begin{lem}
\label{ldc}
Let $(\O,\mu)$ be a measure space. Let $\{\Phi_j\}_{j\ge0}$ be
a sequence of weakly measurable operator-valued functions.
Suppose that $\sup_j\|\Phi_j(\o)\|\in L^1(\mu)$ and
$\Phi(\o)=\lim\limits_{j\to\be}\Phi_j(\o)$ in the strong operator topology for $\mu$-almost all $\o$.
Then
$$
\lim_{j\to\be}\int_\O\Phi_j\,d\mu=\int_\O\Phi\,d\mu
$$
in the strong operator topology.
\end{lem}

{\bf Proof of Theorem \ref{preper0}.}
Clearly, it suffices to prove the result in the case when
$\supp\F f$ is a compact subset of $(0,\be)$. By Lemma \ref{sot},
$$
\lim_{j\to\be}\exp({\rm i}tA_{k,j})=\exp({\rm i}tA_{k})
\quad
\mbox{and}\quad
\lim_{j\to\be}\big((\sm^{-t}f)*q_t\big)(A_{k,j})=\big((\sm^{-t}f)*q_t\big)(A_{k})
$$
in the strong operator topology for all $t>0$. Now the result follows from
Lemma \ref{ldc} and the integral representation for $\dg^{m}f$
as an element of projective tensor product of $L^\be$ spaces, see the discussion in the beginning of this section. $\bl$

\

\section{\bf A formula for operator Taylor polynomials}
\setcounter{equation}{0}
\label{for}

\

In this section we obtain a formula for operator Taylor polynomials in terms multiple operator integrals.
We are going to prove that for $f\in B_{\be1}^m(\R)$, the following formula holds
\bay
\label{of}
{\mathscr T}^{(m)}_{A,K}f=\underbrace{\int\!\!\cdots\!\!\int}_{m+1}(\dg^{m}f)(x_1,\cdots,x_{m+1})
\,dE_{A+K}(x_1)K\,dE_A(x_2)K\cdots K\,dE_A(x_{m+1}).
\ey

We believe that formula \rf{of} is of independent interest though it will not be used to prove the the main results of the paper.

Here $A$ is a (not necessarily bounded) self-adjoint operator and $K$ is a bounded self-adjoint operator.

Recall that
\begin{align}
\label{Tpm}
{\mathscr T}^{(m)}_{A,K}f&\df f(A+K)-f(A)\nonumber\\[.2cm]
&-\frac{d}{dt}f(A+tK)\Big|_{t=0}-
\cdots-\frac1{(m-1)!}\frac{d^{m-1}}{dt^{m-1}}f(A+tK)\Big|_{t=0}.
\end{align}

Actually, it is not true that each term of the right-hand side of \rf{Tpm} exists for functions $f$ in $B_{\be1}^m(\R)$.
However, we explain in this section how to define the right-hand side of \rf{Tpm} for all functions $f$ in $B_{\be1}^m(\R)$.

First, we establish \rf{of} for bounded functions $f$ whose Fourier transform has compact support in $(0,\be)$.

\begin{thm}
\label{cs}
Let $f$ be a function in $L^\be(\R)$ such that the support of $\F f$ is a compact subset of $(0,\be)$. Then
{\em\rf{of}} holds.
\end{thm}

To prove Theorem \ref{cs}, we need the following result.

\begin{lem}
\label{fla}
Let $f$ be a function in $L^\be(\R)$ such that the support of $\F f$ is a compact subset of $(0,\be)$. Then for every $m\ge2$,
\begin{align*}
&\underbrace{\int\!\!\cdots\!\!\int}_{m}\big(\dg^{m-1}f\big)(x_1,\cdots,x_{m})
\,dE_{A+K}(x_1)K\,dE_A(x_2)K\cdots K\,dE_A(x_{m})\\[.2cm]
&-\underbrace{\int\!\!\cdots\!\!\int}_{j}\big(\dg^{m-1}f\big)(x_1,\cdots,x_{m})
\,dE_{A}(x_1)K\,dE_A(x_2)K\cdots K\,dE_A(x_{m})\\[.2cm]
=&\underbrace{\int\!\!\cdots\!\!\int}_{m+1}\big(\dg^{m}f)(x_1,\cdots,x_{m+1}\big)
\,dE_{A+K}(x_1)K\,dE_A(x_2)K\cdots K\,dE_A(x_{m+1}).\\[.2cm]
\end{align*}
\end{lem}

\Pf To simplify the notation, we give the proof in the case $m=2$. The proof in the general case is exactly the same.

Consider first the case when $A$ is bounded. We have
\begin{align*}
&\iint(\dg f)(x,y)\,dE_{A+K}(x)K\,dE_A(y)-\iint(\dg f)(x,y)\,dE_A(x)K\,dE_A(y)\\[.2cm]
=&\iint(\dg f)(x,z)\,dE_{A+K}(x)K\,dE_A(z)-\iint(\dg f)(y,z)\,dE_A(y)K\,dE_A(z)\\[.2cm]
=&\iiint(\dg f)(x,z)\,dE_{A+K}(x)\,dE_A(y)K\,dE_A(z)\\[.2cm]
&-
\iiint(\dg f)(y,z)\,dE_{A+K}(x)\,dE_A(y)K\,dE_A(z)\\[.2cm]
=&\iiint\Big(\big(\dg^2f\big)(x,y,z)\Big)(x-y)\,dE_{A+K}(x)\,dE_A(y)K\,dE_A(z)\\[.2cm]
=&\iiint\Big(\big(\dg^2f\big)(x,y,z)\Big)x\,dE_{A+K}(x)\,dE_A(y)K\,dE_A(z)\\[.2cm]
&-\iiint\Big(\big(\dg^2f\big)(x,y,z)\Big)y\,dE_{A+K}(x)\,dE_A(y)K\,dE_A(z)\\[.2cm]
=&\iiint\big(\dg^2f\big)(x,y,z)\,dE_{A+K}(x)(A+K)\,dE_A(y)K\,dE_A(z)\\[.2cm]
&-\iiint\big(\dg^2f\big)(x,y,z)\,dE_{A+K}(x)A\,dE_A(y)K\,dE_A(z)\\[.2cm]
=&\iiint\big(\dg^2f\big)(x,y,z)\,dE_{A+K}(x)K\,dE_A(y)K\,dE_A(z).
\end{align*}

In the case of unbounded $A$, we put
$$
A_j\df AE_A([-j,j]).
$$
Clearly,
$$
\lim_{j\to\be}\|A_ju-Au\|=0\quad\mbox{for every}\quad u\quad
\mbox{in the domain of}\quad A.
$$
Since each operator $A_j$ is bounded, we have
\begin{align*}
\iint(\dg f)(x,y)&\,dE_{A_j+K}(x)K\,dE_{A_j}(y)-\iint(\dg f)(x,y)\,dE_{A_j}(x)K\,dE_{A_j}(y)\\[.2cm]
=&\iiint\big(\dg^2f\big)(x,y,z)\,dE_{A_j+K}(x)K\,dE_{A_j}(y)K\,dE_{A_j}(z).
\end{align*}
The result follows now from Theorem \ref{preper0}. $\bl$

\medskip

{\bf Proof of Theorem \ref{cs}.} We proceed by induction. It was proved in \cite{Pe1} that identity \rf{of}
holds for $m=1$ and for functions $f$ in the Besov space $B_{\be1}^1(\R)$, and so it holds for functions $f$ satisfying the hypotheses of Theorem \ref{cs}.

To pass from $m-1$ to $m$, we use Lemma \ref{fla} and the formula
\begin{align*}
\frac{d^m}{dt^m}f(A&+tK)\Big|_{t=0}\\[.2cm]
=&
m!\underbrace{\int\!\!\cdots\!\!\int}_{m+1}(\dg^{m}f)(x_1,\cdots,x_{m+1})
\,dE_{A+K}(x_1)K\,dE_A(x_2)K\cdots K\,dE_A(x_{m+1}).
\end{align*}
This formula was proved in \cite{Pe5} for $f\in B_{\be1}^m(\R)$, and so holds for functions satisfying
the hypotheses of Theorem \ref{cs}. $\bl$

We can extend now formula \rf{of} to the case of arbitrary functions $f$ in $B_{\be1}^m(\R)$. As we have already mentioned,
for $f\in B_{\be1}^m(\R)$, each term of the Taylor polynomial $\Tp f$ is not necessarily defined. However,
we can define $\Tp f$ fro $f\in B_{\be1}^m(\R)$ by the following formula:
\bay
\label{tpf}
\Tp f\df\sum_{n=-\be}^\be\Tp f_n,
\ey
where the functions $f_n$ are given by \rf{fn}.  Note that it follows from
\rf{Boke} that the right-hand side of \rf{tpf} converges
absolutely for all functions $f$ in $B_{\be1}^m(\R)$. It is also easy to see that the right-hand side of \rf{tpf} does not depend on the choice of the function $w$ in the definition of Besov spaces, see \S\,2.

\begin{thm}
\label{TfB}
Let $f\in B_{\be1}^m(\R)$. Then {\em\rf{of}} holds, where $\Tp f$ is defined by {\em\rf{tpf}}.
\end{thm}

\Pf By Theorem \ref{of}, \rf{of} holds for functions $f$ in $L^\be(\R)$ such that $\supp\F f$ is
a compact subset of $(0,\be)$. In the same way it can be proved that \rf{of} holds for functions $f$ in
$L^\be(\R)$ such that $\supp\F f$ is a compact subset of $(-\be,0)$. Thus \rf{of} holds for each function $f_n$. The result follows now from Theorem 5.5 of \cite{Pe5}. $\bl$

\

\section{\bf The general result}
\setcounter{equation}{0}
\label{gr}

\

In this section we establish most general spectral formulae for perturbations of class $\bS_m$ (see Theorem \ref{gsl} below). We show in the next section that the trace formula for operator Taylor polynomials is a special case of Theorem \ref{gsl}.

Let $A$ be a self-adjoint operator and let $K$ be a bounded self-adjoint operator.
Put $A_t\df A+tK$ for $t\in\R$. Let $f$ be a bounded continuous function on $\R$ and let $m$ be a positive integer $m$. We can consider the following finite differences:
$$
\big(\D_K^mf\big)(A)\df\sum_{j=0}^m(-1)^{m-j}\left(\begin{matrix}m\\j\end{matrix}\right)f\big(A+jK\big).
$$
It turns out that the finite differences $\big(\D_K^mf\big)(A)$ can be defined for functions $f$ in $B_{\be1}^m(\R)$ by the formula
$$
\big(\D_K^mf\big)(A)\df\sum_{n\in\Z}\big(\D_K^mf_n\big)(A),
$$
where the functions $f_n$ are defined by \rf{fn}. Moreover, it was shown in
\cite{AP1}, Lemma 4.3 that
\begin{align}
\label{kr}
&\big(\D_K^mf\big)(A)\nonumber\\[.2cm]
=&m!\underbrace{\int\!\!\cdots\!\!\int}_{m+1}(\dg^{m}f)(x_1,\cdots,x_{m+1})\,dE_A(x_1)K\,dE_{A+K}(x_2)K\cdots K\,dE_{A+mK}(x_{m+1}).
\end{align}
Strictly speaking, formula \rf{kr} was proved in \cite{AP1} for bounded self-adjoint operators $A$. However, it is easy to see that the approximation procedure used in the proof of Lemma \ref{fla} in this paper also works to extend formula \rf{kr} to the case of unbounded $A$.

Recall that it was shown in \cite{Pe5} that for every function $f$ in $B_{\be1}^m(\R)$, the function $t\mapsto f(A_t)$ has $m$th derivative in the operator norm if we define it by
$$
\frac{d^{m}}{dt^{m}}f(A_t)\df\sum_{n\in\Z}\frac{d^{m}}{dt^{m}}f_n(A_t)
$$
and
\begin{align}
\label{mproi}
&\frac{d^{m}}{dt^{m}}f(A_t)\Big|_{t=s}\nonumber\\[.2cm]
=&m!
\underbrace{\int\cdots\int}_{m+1}(\dg^{m}f)(x_1,\cdots,x_{m+1})
\,dE_{A_s}(x_1)K\,dE_{A_s}(x_2)K\cdots K\,dE_{A_s}(x_{m+1}),
\end{align}
where $E_s$ is the spectral measure of $A_s$.

The following estimate was obtained in \cite{PSS} for functions $f$ satisfying
\rf{PrFo}. We extend it to the class $B_{\be1}^m(\R)$.

\begin{thm}
\label{pss}
Let $f\in B_{\be1}^m(\R)$ and $K\in\bS_m$. Then
\bay
\label{PsS}
\left\|\tr\left(\frac{d^{m}}{dt^{m}}f(A_t)\right)\right\|_{L^\infty}\le
\const\big\|f^{(m)}\big\|_{L^\be}\|K\|^m_{\bS_m}.
\ey
\end{thm}

\Pf Clearly, it suffices to show that
$$
\left|\tr\left(\frac{d^{m}}{dt^{m}}f(A_t)\Big|_{t=0}\right)\right|\le
\const\big\|f^{(m)}\big\|_{L^\be}\|K\|^m_{\bS_m}.
$$

Let us first prove the result in the case when $\supp\F f$ is a compact subset of $(0,\be)$.
Let $\Phi$ be a function in $C^\be(\R)$ such that $\Phi(0)=1$ and $\F\Phi$ is a nonnegative infinitely differentiable function
with a compact support. For $\e>0$, we put $f_\e(x)\df\Phi(\e x)f(x)$.
Then $\supp\F f_\e$ is a compact and
$$
\F f_\e\in L^1(\R)\cap C^\be(\R).
$$
Hence, by
Theorem 2.1 of \cite{PSS},
$$
\left|\tr\left(\frac{d^{m}}{dt^{m}}f_\e(A_t)\Big|_{t=0}\right)\right|\le
\const\big\|f_\e^{(m)}\big\|_{L^\be}\|K\|^m_{\bS_m}.
$$
It is easy to see that
$$
\big\|f_\e^{(m)}\big\|_{L^\be}\le C\big\|f^{(m)}\big\|_{L^\be},
$$
where $C$ depends only on $\Phi$ and $m$.
Moreover, $\supp\F f_\e$ is a compact support of $(0,\be)$ for sufficiently small $\e$
and $\lim_{\e\to0}f_\e=f$ in the space $B_{\be1}^m(\R)$.
It follows from \rf{mproi} and from \rf{Boke} that
$$
\left\|\frac{d^{m}}{dt^{m}}(f-f_\e)(A_t)\right\|_{\bS_1}
\le\const\|f-f_\e\|_{B_{\be1}^m(\R)}\|K\|_{\bS_m}^m.
$$
Hence,
\bay
\label{eps}
\lim_{\e\to0}\tr\left(\frac{d^{m}}{dt^{m}}f_\e(A_t)\Big|_{t=0}\right)=
\tr\left(\frac{d^{m}}{dt^{m}}f(A_t)\Big|_{t=0}\right).
\ey
%

If $\supp\F f$ is a compact support of $(-\be,0)$, the proof of \rf{eps} is the same.

It follows from \rf{eps} that for $f$ in $B_{\be1}^m(\R)$,
$$
\left|\tr\left(\frac{d^{m}}{dt^{m}}f_n(A_t)\Big|_{t=0}\right)\right|\le
\const\big\|f_n^{(m)}\big\|_{L^\be}\|K\|^m_{\bS_m}.
$$
where the functions $f_n$ are defined by \rf{fn}. This implies \rf{PsS}. $\bl$


\begin{thm}
\label{dfye}
Let $m$ be a positive integer and let $f\in B_{\be1}^m(\R)$. Then
\bay
\label{mth}
\frac{d^{m}}{dt^{m}}f(A_t)\Big|_{t=s}=\lim_{h\to0}h^{-m}\big(\D_{hK}^mf\big)(A_s).
\ey
\end{thm}

\Pf The result follows from \rf{mproi}, \rf{kr} and the equality
\begin{align*}
\lim_{h\to\0}&\underbrace{\int\cdots\int}_{m+1}(\dg^{m}f)(x_1,\cdots,x_{m+1})\,dE_A(x_1)K\,dE_{A+hK}(x_2)K\cdots K\,dE_{A+mhK}(x_{m+1})\\[.2cm]
=&\underbrace{\int\cdots\int}_{m+1}(\dg^{m}f)(x_1,\cdots,x_{m+1})\,dE_A(x_1)K\,dE_{A}(x_2)K\cdots K\,dE_{A}(x_{m+1}),
\end{align*}
which is a special case of Theorem \ref{preper}. $\bl$

Consider  now the case when $K\in\bS_m$.

\begin{thm}
\label{mbesov}
Let $f\in B^m_{\be1}(\R)$ and $K\in\bS_m$.
Then the limit
\bay
\label{pred}
\lim_{h\to0}h^{-m}\big(\D_{hK}^mf\big)(A_s)
\ey
exists in the norm of $\bS_1$ for every $s$. Moreover,
the function
\bay
\label
{kont}
s\mapsto \frac{d^{m}}{dt^{m}}f(A_t)\Big|_{t=s}
\ey
is a continuous $\bS_1$-valued function and
\bay
\label{nervo}
\left\|\frac{d^{m}}{dt^{m}}f(A_t)\right\|_{L^\be(\bS_1)}\le \const\|f\|_{B^m_{\be1}}\|K\|^m_{\bS_m}.
\ey
\end{thm}

\Pf The fact that the limit \rf{pred} exists in the norm of $\bS_1$ is an immediate consequence of Theorem \ref{preper}. The continuity of the function
\rf{kont} also immediately follows from Theorem \ref{preper}. Finally, inequality \rf{nervo} follows from \rf{mproi} and inequality \rf{Besov}.
$\bl$

Let $C_0(\R)$ denote the space of all functions $f\in C(\R)$ such that
$\lim_{|t|\to\be}f(t)=0$.
Denote by $\M(\R)$ the space of all finite Borel measures on $\R$.
We identify in a natural way the space $\M(\R)$ with the space $(C_0(\R))^*$.

\begin{thm}
\label{mery}
Let $m$ be a positive interger and
let $A$ and $K$ be self-adjoint operators such that $K\in\bS_m$.
Then for each $t\in\R$ there exists a unique measure $\nu_t\in\M(\R)$
such that
\bay
\label{nut}
\|\nu_t\|\le \const\|K\|_{\bS_m}\quad\mbox{and}
\quad\tr\left( \frac{d^{m}}{dt^{m}}f(A_t)\right)
=\int_\R f^{(m)}\,d\nu_t
\ey
for every $f\in B^m_{\be1}(\R)$. Moreover, the map $t\mapsto\nu_t$ is a continuous map from $\R$ to
$\M(\R)$ equipped with the weak-$*$ topology $\s(\M(\R),C_0(\R))$.
\end{thm}

\Pf Note that the set $\frak X\df B^0_{\be1}(\R)\cap C_0(\R)$
is dense in $C_0(\R)$.
Theorem \ref{pss} implies that there exists a unique measure $\nu_t\in\M(\R)$
such that \rf{nut} holds for every
$f\in\frak X$. Theorem \ref{mbesov} implies that the mapping $t\mapsto\nu_t$
is continuous.

It remains to establish the equality in \rf{nut} for all
$f\in B^m_{\be1}(\R)$. Clearly, it suffices to verify that this equality holds on a dense
subset of $B^m_{\be1}(\R)$. Let $f$ be an entire function of exponential type $\s$
that is bounded on $\R$.
Put $\f(t)\df t^{-1}\sin t$. It is easy to see that
$$
\lim_{|t|\to\be}\big(\f(\e t)f(t)\big)^{(m)}=0,
$$
$$
\lim_{\e\to0}\big(\f(\e t)f(t)\big)^{(m)}=\big(f(t)\big)^{(m)}(t),
$$
and
$$
\sup_{\e,t}\left|\big(\f(\e t)f(t)\big)^{(m)}\right|<\be.
$$


To complete the proof, it suffices to verify that
$$
\lim_{\e\to0}\tr\left(\frac{d^{m}}{dt^{m}}\big(\f(A_{\e t})\,f(A_t)\big)\right)=
\tr\left(\frac{d^{m}}{dt^{m}}f(A_t)\right).
$$
We have
$$
\frac{d^{m}}{dt^{m}}\big(\f(A_{\e t})\,f(A_t)\big)
=\sum_{j=0}^m {m\choose j}\,\frac{d^{j}}{dt^{j}}\f(A_{\e t})\,\,\frac{d^{m-j}}{dt^{m-j}}f(A_t).
$$
Inequality (5.3) in \cite{AP2} implies the following inequalities:
$$
\left\|\frac{d^{j}}{dt^{j}}\f(A_{\e t})\right\|_{\bS_{\frac mj}}\le C_m\e^{j}\|K\|_{S_m}^j
\quad\text {and}\quad\left\|\frac{d^{m-j}}{dt^{m-j}}f(A_t)\right\|_{\bS_{\frac m{m-j}}}\le C_m\s^{m-j}\|K\|_{S_m}^{m-j}.
$$
Hence,
$$
\lim_{\e\to0}\left\|\frac{d^{j}}{dt^{j}}\f(A_{\e t})\,\,\frac{d^{m-j}}{dt^{m-j}}f(A_t)\right\|_{\bS_1}=0,\quad j\ge1.
$$
It remains to prove that
$$
\lim_{\e\to0}\tr\left(\f(A_{\e t})\,\frac{d^{m}}{dt^{m}}f(A_t)\right)=\tr\left(\frac{d^{m}}{dt^{m}}f(A_t)\right),
$$
which is clear, because $\lim\limits_{\e\to0}\f(A_{\e t})=I$ in the strong operator
topology. $\bl$

For $\mu\in\M(\R)$, we denote by $\mu^{(m)}$ the $m$-th derivative of $\mu$ in
the sense of distributions.
Thus we put
$$
\left\langle f(A_t),\mu^{(m)}\right\rangle\df(-1)^m\int_\R \frac{d^{m}}{dt^{m}}f(A_t)\,d\mu
\quad\mbox{for}\quad f\in B^m_{\be1}(\R).
$$
Let
\bay
\label{AK}
\mu_{\{A,K\}}\df(-1)^m\int_\R \nu_t\,d\mu(t),
\ey
the integral exists in $\M(\R)$ in the topology $\s(\M(\R),C_0(\R))$, because the function $t\mapsto\nu_t$
is weak star continuous. 

We prove in the next section that if $\mu$ is an absolutely continuous measure (with respect to Lebesgue measure), then $\mu_{\{A,K\}}$ must also be absolutely continuous.

The following result is a most general trace formula for perturbations of class $\bS_m$.

\begin{thm}
\label{gsl}
Let $m$ be a positive integer and
let $A$ and $K$ be self-adjoint operators such that $K\in\bS_m$. Suppose that $\mu\in\M(\R)$. Then for every $f\in B^m_{\be1}(\R)$,
\bay
\label{msled}
\tr\left\langle f(A_t),\mu^{(m)}\right\rangle=\int_\R f^{(m)}\,d\mu_{\{A,K\}},
\ey
where $\mu_{\{A,K\}}$ is defined by
{\rm\rf{AK}}.
\end{thm}

\Pf We have
\begin{align*}
\tr\left\langle f(A_t),\mu^{(m)}\right\rangle&=(-1)^m\int_\R \tr\left(\frac{d^{m}}{dt^{m}}f(A_t)\right)\,d\mu\\[.2cm]
&=(-1)^m\int_\R\left(\int_\R f^{(m)}\,d\nu_t\right)\,d\mu=\int_\R f^{(m)}\,d\mu_{\{A,K\}}.\quad\bl
\end{align*}

\

\section{\bf Trace formulae for operator Taylor polynomials and other special cases}
\setcounter{equation}{0}
\label{sle}

\

In this section we show that trace formula \rf{PSS} is a special case of our general formula \rf{msled} which allows us to improve the main results of \cite{PSS} and extend trace formula \rf{PSS} to the case of functions in the Besov space $B_{\be1}^m(\R)$. We also consider several other interesting special cases of formula \rf{msled}.

\begin{thm}
\label{Tayp}
Let $m$ be a positive integer and
let $A$ and $K$ be self-adjoint operators such that $K\in\bS_m$.
Consider the the
absolutely continuous measure $\mu$ defined by
\bay
\label{mera}
d\mu(t)\df\frac{(-1)^m}{(m-1)!}(1-t)^{m-1}\chi_{[0,1]}(t)\,dt.
\ey
then the formula
$$
\trace\left({\mathscr T}^{(m)}_{A,K}f\right)=\int_\R f^{(m)}(x)\,d\mu_{\{A,K\}}(x)
$$
holds for every function $f$ in the Besov space $B_{\be1}^m(\R)$.
\end{thm}

\Pf
It is easy to see that
$$
\mu^{(m)}=\d_1-\sum_{j=0}^{m-1}\frac{(-1)^j}{j!}\d_0^{(j)},
$$
where $\d_a$ is the unit point mass at $a$.
Indeed, by Taylor's formula,
$$
\f(1)=\sum_{j=0}^{m-1}\frac{\f^{(j)}(0)}{j!}+\int_0^1\frac{(1-t)^{m-1}}{(m-1)!}\f^{(m)}(t)\,dt
$$
for every $\f\in C^m(\R)$. It follows that
$$
{\mathscr T}^{(m)}_{A,K}f=\left\langle f(A_t),\mu^{(m)}\right\rangle.
$$
The result follows now from Theorem \ref{gsl}. $\bl$

\begin{cor}
\label{flasl}
Under the hypotheses of Theorem {\em\ref{Tayp}}, the spectral shift function $\eta_m$ of order $m$ is the Radon--Nykodym derivative of $\mu_{\{A,K\}}$ with respect to Lebesgue measure,
$$
\frac{d\mu_{\{A,K\}}}{d\m}=\eta_m,
$$
and trace formula {\em\rf{PSS}} holds for every $f\in B_{\be1}^m(\R)$.
\end{cor}

\Pf It follows from Theorem \ref{mera} and from trace formula \rf{PSS} that
$$
\int_\R f^{(m)}(x)\,\mu_{\{A,K\}}=\int_\R f^{(m)}(x)\eta_m(x)\,dx
$$
for every infinitely smooth functions $f$ with compact support. This implies the result. $\bl$

Now we are in a position to prove that $\mu_{\{A,K\}}$ is absolutely continuous whenever $\mu$ is.

\begin{thm}
\label{abne}
Suppose that $\mu$ is an absolutely continuous measure (with respect to Lebesgue measure). Then $\mu_{\{A,K\}}$ is also absolutely continuous.
\end{thm}

\Pf For $g\in L^1(\m)$, we use the notation
$$
g_{\{A,K\}}=\mu_{\{A,K\}},
$$
where $\mu$ is the absolutely continuous measure defined by $d\mu=g\,d\m$.

Denote by $\cL$ the set of all $g\in L^1(\m)$ such that the measure
$g_{\{A,K\}}$ is absolutely continuous for all self-adjoint operators $A$ and $K$ with $K\in\bS_m$. Clearly, $\cL$ is a closed translation and dilation invariant
subspace of $L^1(\m)$. By Theorem \ref{Tayp} and Corollary \ref{flasl}, the function
$t\mapsto (1-t)^{m-1}\chi_{[0,1]}(t)$ belongs to $\cL$. 
Since $\int_0^1(1-t)^{m-1}\,dt\not=0$, it follows that $\cL=L^1(\m)$. $\bl$

The following results are interesting special cases of our general trace formula \rf{msled}.

\begin{thm}
\label{fsl2}
Let $m$ be a positive integer.
Suppose that $A$ is a self-adjoint operator and $K$ is a self-adjoint operator of class $\bS_m$.
Then there exists a function $\vk_m\in L^1(\R)$ such that
$$
\trace \big(\D_K^mf\big)(A)=\int_\R f^{(m)}(x)\vk_m(x)\,dx
$$
for every $f\in B_{\be1}^m(\R)$.
\end{thm}

We need the following lemma.
%
%
%

\begin{lem}
\label{l1}
Let $m$ be a positive integer and
let $\l$ be a measure on  a finite interval $[a,b]$. Suppose that $\l$  is orthogonal
to the polynomials of degree less than $m$. Then there exists a
function $g\in L^\be(\R)$ such that $g^{(m)}=\l$ in the sense of distributions, $\supp g\subset[a,b]$
and $\|g\|_{L^1}\le (m!)^{-1}(b-a)^m\|\l\|_{\M(\R)}$.
\end{lem}

\Pf Without loss of generality we may assume that $a=0$. Put
$$
g(x)=\int_0^x\frac{(x-t)^{m-1}}{(m-1)!}\,d\l(t).
$$
It remains to observe that $g(x)=0$ for $x\not\in[0,b]$ and
$$
|g(x)|\le\frac{x^{m-1}}{(m-1)!}\|\l\|_{\M(\R)}
$$
for $x\in(0,b)$. $\bl$

\medskip

{\bf Proof of Theorem \ref{fsl2}.} Put
$$
\l\df \sum_{j=0}^m(-1)^{m-j}{m\choose j}\d_j.
$$
By Lemma \ref{l1}, there exists a function $g\in L^\be(\R)$ such that $\supp g\subset[0,m]$
and $\mu=g^{(m)}$
It is easy to see that
$$
\big(\D_K^mf\big)(A)=\int_\R f(A_t)\,d\l(t).
$$
Hence, by Theorem \ref{gsl}
we have
$$
\tr\big(\big(\D_K^mf\big)(A)\big)=\int_\R f^{(m)}\,d\mu_{\{A,K\}},
$$
where $d\mu(t)=g(t)\,dt$. Moreover, the measure $\mu_{\{A,K\}}$
is absolutely continuous, see \S\,\ref{gr}. $\bl$

Now we state the following generalization of Theorem \ref{fsl2}.

\begin{thm}
\label{fsl3}
Let $m$ be a positive integer and let $A$ and $K$ be self-adjoint operators such that $K\in\bS_m$.
Suppose that the measure $\sum_{j=0}^N\l_j\d_{t_j}$
is orthogonal to the polynomials of degree less than $m$. Then there exists a function $\f$ in $L^1(\R)$ such that
$$
\trace \left(\sum_{j=0}^N\l_jf(A_{t_j})\right)=\int_\R f^{(m)}(x)\f(x)\,dx
$$
for every $f\in B_{\be1}^m(\R)$.
\end{thm}

The proof of Theorem \ref{fsl3} is similar to that of Theorem \ref{fsl2}
and we omit it.

Let us now generalize Theorem \ref{fsl3}.

\begin{thm}
\label{fsl4}
Let $m$ and $k$ be nonnegative integers such that $0\le k<m$. Suppose that the measure $\sum_{j=0}^N\l_j\d_{t_j}$
is orthogonal to the polynomials of degree less than $m-k$. Let
$A$ be a self-adjoint operator and $K$ be a self-adjoint operator of class $\bS_m$.
Then there exists a function $\psi\in L^1(\R)$ such that
$$
\trace \left(\sum_{j=0}^N\l_j\frac{d^{m_0}}{dt^{m_0}}
f(A_t)\Big|_{t=t_j}\right)=\int_\R f^{(m)}(x)\psi(x)\,dx
$$
for every $f\in B_{\be1}^m(\R)$.
\end{thm}

\Pf
By Lemma \ref{l1}, there exists a function $g\in L^\be(\R)$ with a compact support such that
and $g^{(m-m_0)}=\sum_{j=0}^N\l_j\d_{t_j}$. Then $g^{(m)}=\sum_{j=0}^N\l_j\d_{t_j}^{(m_0)}$.
By Theorem \ref{gsl},
we have
$$
\sum_{j=0}^N\l_j\frac{d^{m_0}}{dt^{m_0}}f(A_t)\Big|_{t=t_j}=\int_\R f^{(m)}\,d\mu_{\{A,K\}},
$$
where $\mu_{\{A,K\}}$ is defined by \rf{AK} with $d\mu(t)=g(t)\,dt$. Clearly, the measure $\mu_{\{A,K\}}$
is absolutely continuous with respect to Lebesgue measure. $\bl$

We have considered several special cases our general trace formula \rf{gsl}. In all those cases the measure $\mu_{\{A,K\}}$ is absolutely continuous with respect to Lebesgue measure. The following example shows that this is not always the case.

\medskip

{\bf Example.} Suppose that $A$ and $K$ are commuting
self-adjoint operators with discrete spectra, i.e.,
$A=\sum_j\a_j(\cdot,e_j)e_j$ and $K=\sum_j\l_j(\cdot,e_j)e_j$,
where $\{e_j\}_j$ is
an orthonormal basis  and the $\a_j$ and $\mu_j$ are real numbers.
By Theorem \ref{mery},
$$
\tr\left( \frac{d^{m}}{dt^{m}}f(A_t)\right)
=\int_\R f^{(m)}\,d\nu_t
$$
for a unique measure $\nu_t$ on $\R$.
Note that this is a special case of trace formula \rf{gsl}.
It is easy to verify that
$$
\tr\left( \frac{d^{m}}{dt^{m}}f(A_t)\right)=
\frac{d^{m}}{dt^{m}}\tr\big(f(A_t)\big)=
\sum_j\l_j^mf^{(m)}(\a_j+t\l_j)
$$
for every $f\in C^m(\R)$ with $f^{(m)}\in L^\be(\R)$.
Thus
$$
\nu_t=\sum_j\l_j^m\d_{\a_j+t\l_j}.
$$
Hence, $\nu_t$ is a discrete measure.

This example shows that in trace formula \rf{msled} the measure $\mu_{\{A,K\}}$ does not have to be absolutely continuous in general.

\

\

\noindent
\begin{tabular}{p{9cm}p{15cm}}
A.B. Aleksandrov & V.V. Peller \\
St-Petersburg Branch & Department of Mathematics \\
Steklov Institute of Mathematics  & Michigan State University \\
Fontanka 27, 191023 St-Petersburg & East Lansing, Michigan 48824\\
Russia&USA
\end{tabular}

\end{document}